\documentclass[12pt]{amsart}

\usepackage{amsmath,amsfonts,amssymb,amscd,verbatim,delarray,verbatim}
\DeclareFontFamily{OT1}{wncyr}{\hyphenchar\font45}
\DeclareFontShape{OT1}{wncyr}{m}{n}{%
   <5> <6> <7> <8> <9> gen * wncyr
   <10> <10.95> <12> <14.4> <17.28> <20.74>  <24.88>wncyr10}{}
\DeclareFontShape{OT1}{wncyr}{m}{it}{%
   <5> <6> <7> <8> <9> gen * wncyi
   <10> <10.95> <12> <14.4> <17.28> <20.74> <24.88> wncyi10}{}
\DeclareFontShape{OT1}{wncyr}{m}{sc}{%
   <5> <6> <7> <8> <9> <10> <10.95> <12> <14.4>
   <17.28> <20.74> <24.88>wncysc10}{}
\DeclareFontShape{OT1}{wncyr}{b}{n}{%
   <5> <6> <7> <8> <9> gen * wncyb
   <10> <10.95> <12> <14.4> <17.28> <20.74> <24.88>wncyb10}{}
\input cyracc.def
\def\rus{\usefont{OT1}{wncyr}{m}{n}\cyracc\fontsize{9}{11pt}\selectfont}

\DeclareMathSizes{9}{9}{7}{5} 

\newcommand{\C}{{\mathbb C}}
\newcommand{\F}{{\mathbb F}}
\newcommand{\G}{{\mathbb G}}

\newcommand{\Q}{{\mathbb Q}}
\newcommand{\Z}{{\mathbb Z}}

\newcommand{\Sha}{\mbox{\rus{\fontsize{11}{11pt}\selectfont{SH}}}}

\DeclareMathOperator{\Aut}{Aut}

\DeclareMathOperator{\Pic}{Pic\,}
\DeclareMathOperator{\Br}{Br\,}

\begin{document}
\numberwithin{equation}{section}

\newtheorem{theorem}{Theorem}[section]
\newtheorem{lemma}[theorem]{Lemma}

\newtheorem{prop}[theorem]{Proposition}
\newtheorem{proposition}[theorem]{Proposition}
\newtheorem{corollary}[theorem]{Corollary}
\newtheorem{corol}[theorem]{Corollary}
\newtheorem{conj}[theorem]{Conjecture}
\newtheorem{sublemma}[theorem]{Sublemma}

\theoremstyle{definition}
\newtheorem{defn}[theorem]{Definition}
\newtheorem{example}[theorem]{Example}
\newtheorem{examples}[theorem]{Examples}
\newtheorem{remarks}[theorem]{Remarks}
\newtheorem{remark}[theorem]{Remark}
\newtheorem{algorithm}[theorem]{Algorithm}
\newtheorem{question}[theorem]{Question}
\newtheorem{problem}[theorem]{Problem}
\newtheorem{subsec}[theorem]{}
\newtheorem{clai}[theorem]{Claim}

\def\toeq{{\stackrel{\sim}{\longrightarrow}}}
\def\into{{\hookrightarrow}}


\def\alp{{\alpha}}  \def\bet{{\beta}} \def\gam{{\gamma}}
 \def\del{{\delta}}
\def\eps{{\varepsilon}}
\def\kap{{\kappa}}                   \def\Chi{\text{X}}
\def\lam{{\lambda}}
 \def\sig{{\sigma}}  \def\vphi{{\varphi}} \def\om{{\omega}}
\def\Gam{{\Gamma}}   \def\Del{{\Delta}}
\def\Sig{{\Sigma}}   \def\Om{{\Omega}}
\def\ups{{\upsilon}}


\def\F{{\mathbb{F}}}
\def\BF{{\mathbb{F}}}
\def\BN{{\mathbb{N}}}
\def\Q{{\mathbb{Q}}}
\def\Ql{{\overline{\Q }_{\ell }}}
\def\CC{{\mathbb{C}}}
\def\R{{\mathbb R}}
\def\V{{\mathbf V}}
\def\D{{\mathbf D}}
\def\BZ{{\mathbb Z}}
\def\K{{\mathbf K}}
\def\XX{\mathbf{X}^*}
\def\xx{\mathbf{X}_*}

\def\AA{\Bbb A}
\def\BA{\mathbb A}
\def\HH{\mathbb H}
\def\PP{\Bbb P}

\def\Gm{{{\mathbb G}_{\textrm{m}}}}
\def\Gmk{{{\mathbb G}_{\textrm m,k}}}
\def\GmL{{\mathbb G_{{\textrm m},L}}}
\def\Ga{{{\mathbb G}_a}}

\def\Fb{{\overline{\F }}}
\def\Kb{{\overline K}}
\def\Yb{{\overline Y}}
\def\Xb{{\overline X}}
\def\Tb{{\overline T}}
\def\Bb{{\overline B}}
\def\Gb{{\bar{G}}}
\def\Ub{{\overline U}}
\def\Vb{{\overline V}}
\def\Hb{{\bar{H}}}
\def\kb{{\bar{k}}}

\def\Th{{\hat T}}
\def\Bh{{\hat B}}
\def\Gh{{\hat G}}

\def\cF{{\mathfrak{F}}}
\def\cO{{\mathcal O}}
\def\cU{{\mathcal U}}

\def\Lt{{\widetilde L}}
\def\Gt{{\widetilde G}}

\def\gg{{\mathfrak g}}
\def\hh{{\mathfrak h}}
\def\lie{\mathfrak a}

\def\XX{\mathfrak X}
\def\RR{\mathfrak R}
\def\NN{\mathfrak N}

\def\minus{^{-1}}

\def\GL{\textrm{GL}}            \def\Stab{\textrm{Stab}}
\def\Gal{\textrm{Gal}}          
\def\Lie{\textrm{Lie\,}}        \def\Ext{\textrm{Ext}}
\def\PSL{\textrm{PSL}}          \def\SL{\textrm{SL}}
\def\loc{\textrm{loc}}
\def\coker{\textrm{coker\,}}    \def\Hom{\textrm{Hom}}
\def\im{\textrm{im\,}}           \def\int{\textrm{int}}
\def\inv{\textrm{inv}}           \def\can{\textrm{can}}
\def\id{\textrm{id}}              \def\Char{\textrm{char}}
\def\Cl{\textrm{Cl}}
\def\Sz{\textrm{Sz}}
\def\ad{\textrm{ad\,}}
\def\SU{\textrm{SU}}
\def\Sp{\textrm{Sp}}
\def\PSL{\textrm{PSL}}
\def\PSU{\textrm{PSU}}
\def\rk{\textrm{rk}}
\def\PGL{\textrm{PGL}}
\def\Ker{\textrm{Ker}}
\def\Ob{\textrm{Ob}}
\def\Var{\textrm{Var}}
\def\poSet{\textrm{poSet}}
\def\Al{\textrm{Al}}
\def\Int{\textrm{Int}}
\def\Mod{\textrm{Mod}}
\def\Smg{\textrm{Smg}}
\def\ISmg{\textrm{ISmg}}
\def\Ass{\textrm{Ass}}
\def\Grp{\textrm{Grp}}
\def\Com{\textrm{Com}}
\def\rank{\textrm{rank}}

\def\char{\textrm{char}}

\newcommand{\Or}{\operatorname{O}}

\def\tors{_\def{\textrm{tors}}}      \def\tor{^{\textrm{tor}}}
\def\red{^{\textrm{red}}}         \def\nt{^{\textrm{ssu}}}

\def\sss{^{\textrm{ss}}}          \def\uu{^{\textrm{u}}}
\def\mm{^{\textrm{m}}}
\def\tm{^\times}                  \def\mult{^{\textrm{mult}}}

\def\uss{^{\textrm{ssu}}}         \def\ssu{^{\textrm{ssu}}}
\def\comp{_{\textrm{c}}}
\def\ab{_{\textrm{ab}}}

\def\et{_{\textrm{\'et}}}
\def\nr{_{\textrm{nr}}}

\def\nil{_{\textrm{nil}}}
\def\sol{_{\textrm{sol}}}
\def\End{\textrm{End\,}}

\def\til{\;\widetilde{}\;}

\def\min{{}^{-1}}



\title[Algebraic tori --- thirty years after]
{{\bf Algebraic tori --- thirty years after}}

\author[Boris Kunyavskii] {Boris Kunyavski\u\i}

\address{Department of
Mathematics, Bar-Ilan University, 52900 Ramat Gan, ISRAEL}
\email{kunyav@macs.biu.ac.il}

\dedicatory{To my teacher Valentin Evgenyevich Voskresenski\u\i , with gratitude
and admiration}


\begin{abstract}
This is an expanded version of my talk given at the International
Conference ``Algebra and Number Theory'' dedicated to the 80th
anniversary of V.~E.~Voskresenskii, which was held  at the Samara
State University in May 2007. The goal is to give an overview of
results of V.~E.~Voskresenskii  on arithmetic and birational
properties of algebraic tori which culminated in his monograph
``Algebraic Tori'' published in Russian 30 years ago. I shall try to
put these results and ideas into somehow broader context and also to
give a brief digest of the relevant activity related to the period
after the English version of the monograph ``Algebraic Groups and
Their Birational Invariants'' appeared.
\end{abstract}

\maketitle

\section*{} \label{intro}

\section{Rationality and nonrationality problems} \label{sec:rat}

A classical problem, going back to Pythagorean triples, of describing
the set of solutions of a given system of polynomial equations by
rational functions in a certain number of parameters ({\it rationality
problem}) has been an
attraction for many generations. Although a lot of various techniques
have been used, one can notice that after all, to establish rationality,
one usually has to exhibit some explicit parameterization such as that
obtained by stereographic projection in the Pythagoras problem. The situation
is drastically different if one wants to establish non-existence of such a
parameterization ({\it nonrationality problem}): here one usually has to use
some known (or even invent some new) {\it birational invariant} allowing one
to detect nonrationality by comparing its value for the object under consideration
with some ``standard'' one known to be zero; if the computation gives a nonzero
value, we are done. Evidently, to be useful, such an invariant must be (relatively)
easily computable. Note that the above mentioned subdivision to rationality and
nonrationality problems is far from being absolute: given a class of objects,
an ultimate goal could be to introduce some computable birational invariant giving
{\it necessary and sufficient} conditions for rationality. Such a task may be not
so hopeless, and some examples will be given below.

In this section, we discuss several rationality and nonrationality problems related
to {\it {algebraic tori}}. Since this is the main object of our consideration,
for the reader's convenience we shall recall the definition.

\begin{defn} \label{def:tor}
Let $k$ be a field. An algebraic $k$-torus $T$ is an algebraic $k$-group such that
over a (fixed) separable closure $\bar k$ of $k$ it becomes isomorphic to a direct
product of $d$ copies of the multiplicative group:
$$
T\times _k\bar k\cong \G_{\text{\rm{m}},\bar k}^d
$$
(here $d$ is the dimension of $T$).
\end{defn}

We shall repeatedly use the duality of categories of
algebraic $k$-tori and finite-dimensional torsion-free $\Z$-modules
viewed together with the action of the Galois group $\gg :=\Gal (\kb /k)$ which is given
by associating to a torus $T$ its $\gg$-module of characters
$\Th :=\Hom (T\times \kb, \G_{\text{\rm{m}},\kb })$. Together with the fact
that $T$ splits over a finite extension of $k$ (we shall denote by $L$ the smallest
of such extensions and call it the {\it minimal splitting field} of $T$), this allows
us to reduce many problems to considering (conjugacy classes of) finite subgroups
of $GL(d,\Z )$ corresponding to $\Th$.

\subsection{Tori of small dimension}

There are only two subgroups in $GL(1,\Z )$: $\{(1)\}$ and $\{(1),(-1)\}$, both corresponding
to $k$-rational tori: $\G_{{\text{\rm{m}}},k}$ and $R_{L/k}\G_{{\text{\rm{m}}},L}/\G_{{\text{\rm{m}}},k}$,
respectively (here $L$ is a separable quadratic extension of $k$ and $R_{L/k}$ stands for the Weil functor of
restriction of scalars).

For $d=2$ the situation was unclear until the breakthrough obtained by Voskresenski\u\i \ \cite{Vo67} who proved
\begin{theorem} \label{th:two}
All two-dimensional tori are $k$-rational.
\end{theorem}

For $d=3$ there exist nonrational tori, see \cite{Ku87} for birational classification.

For $d=4$ there is no classification, some birational invariants were computed in \cite{Po}.

\begin{remark} \label{rem-class-free}
The original proof of rationality of the two-dimensional tori in \cite{Vo67} is based on the
classification of finite subgroups in $GL(2,\Z )$ and case-by-case analysis. In the monographs \cite{Vo-rus} and
\cite{Vo-eng} a simplified proof is given, though using the classification of {\it maximal} finite subgroups
in $GL(2,\Z )$. Merkurjev posed a question about the existence of a classification-free proof. (Note that one of
his recent results on 3-dimensional tori \cite{Me} (see Section \ref{sec:R} below) was obtained without
referring to the birational classification given in \cite{Ku87}.) During discussions in Samara, Iskovskikh and Prokhorov
showed me a proof only relying on an ``easy'' part of the classification theorem for rational surfaces.
\end{remark}

\begin{remark} \label{rem:cremona}
One should also note a recent application of Theorem \ref{th:two} to the problem of classification of elements
of prime order in the Cremona group $Cr(2,\Q )$ of birational transformations of plane (see \cite[6.9]{Se}). The problem
was recently solved by Dolgachev and Iskovskikh \cite{DI}. See the above cited papers for more details.
\end{remark}

\subsection{Invariants arising from resolutions} \label{sec:res}

Starting from the pioneering works of Swan and Voskresenski\u\i \ on Noether's problem (see Section \ref{sec:N}),
it became clear that certain resolutions of the Galois module $\Th$ play an important role in understanding
birational properties of $T$. These ideas were further developed by Lenstra, Endo and Miyata,
Colliot-Th\'el\`ene and Sansuc (see \cite{Vo-rus}, \cite{Vo-eng} for an account of that period); they were pursued
through several decades, and some far-reaching generalizations were obtained in more recent works (see the remarks at the end of this
and the next sections).

For further explanations we need to recall some definitions. Further on ``module'' means a finitely generated
$\Z$-free $\gg$-module.

\begin{defn} \label{def:fl}
We say that $M$ is a {\it permutation} module if it has a $\Z$-basis permuted by $\gg$.
We say that modules $M_1$ and $M_2$ are similar if there exist permutation modules $S_1$ and $S_2$
such that $M_1\oplus S_1\cong M_2\oplus S_2$. We denote the similarity class of $M$ by $[M]$.
We say that $M$ is a {\it coflasque} module if $H^1(\hh , M)=0$ for all closed subgroups $\hh$
of $\gg$. We say that $M$ is a {\it flasque} module if its dual module $M^{\circ}:=\Hom (M,\Z )$ is
coflasque.
\end{defn}

The following fact was first established in the case when $k$ is a field of characteristic zero by
Voskresenski\u\i \ \cite{Vo69} by a geometric construction (see below) and then in \cite{EM}, \cite{CTS77}
in a purely algebraic way (the uniqueness of $[F]$ was established independently by all above authors).

\begin{theorem} \label{th:fl}
Any module $M$ admits a resolution of the form
\begin{equation}
0\to M\to S \to F \to 0,
\label{eq:fl}
\end{equation}
where $S$ is a permutation module and $F$ is a flasque module. The similarity
class $[F]$ is determined uniquely.
\end{theorem}
We denote $[F]$ by $p(M)$. We call (\ref{eq:fl}) the flasque resolution of $M$.
If $T$ is a $k$-torus with character module $M$, the sequence of tori dual to
$(\ref{eq:fl})$ is called the flasque resolution of $T$.

Seeming strange at the first glance, the notion of flasque module (and the
corresponding flasque torus) turned out to be very useful. Its
meaning is clear from the following theorem due to Voskresenski\u\i \ \cite[4.60]{Vo-rus}:

\begin{theorem} \label{th:stable}
If tori $T_1$ and $T_2$ are birationally equivalent, then $p(\hat T_1)=p(\hat T_2)$.
Conversely, if  $p(\hat T_1)=p(\hat T_2)$, then $T_1$ and $T_2$ are {\emph stably equivalent},
i.e. $T_1\times  \G_{\text{\rm{m}},k}^{d_1}$ is birationally equivalent to $T_2\times \G_{\text{\rm{m}},k}^{d_2}$
for some integers $d_1$, $d_2$.
\end{theorem}

Theorem \ref{th:stable} provides a birational invariant of a torus which can be computed in a
purely algebraic way. Moreover, it yields other invariants of cohomological nature which are even easier to compute.
The most important among them is $H^1(\gg ,F)$. One should note that it is well-defined in light of Shapiro's lemma
(because $H^1(\gg ,S)=0$ for any permutation module $S$) and, in view of the inflation-restriction
sequence, can be computed at a finite level: if $L/k$ is the minimal splitting field of $T$ and $\Gamma=\Gal (L/k)$,
we have $H^1(\gg,F)=H^1(\Gamma ,F)$. The finite abelian group $H^1(\Gamma ,F)$ is extremely important for
birational geometry and arithmetic, see below.

\begin{remark} \label{rem:fl-res}
Recently Colliot-Th\'el\`ene \cite{CT07} discovered a beautiful generalization of the flasque resolution
in a much more general context, namely, for an arbitrary connected linear algebraic group $G$.
\end{remark}

\subsection{Geometric interpretation of the flasque resolution}

As mentioned above, the flasque resolution (\ref{eq:fl}) was originally constructed
in a geometric way. Namely, assuming $k$ is of characteristic zero, one can use
Hironaka's resolution of singularities to embed a given $k$-torus $T$ into a smooth
complete $k$-variety $V$ as an open subset and consider the exact sequence of $\gg$-modules
\begin{equation}
0\to \hat T \to S \to \Pic \Vb \to 0,
\label{eq:Pic}
\end{equation}
where $\Vb =V\times_k\kb$, $\Pic \Vb$ is the Picard module, and $S$ is
a permutation module (it is generated by the components of the divisors of $\Vb$ whose
support is outside $\overline T$). With another choice of an embedding $T\to V'$ we are led
to an isomorphism $\Pic \Vb \oplus S_1\cong  \Pic \Vb '\oplus S_2$, so the similarity
class $[\Pic \Vb ]$ is well defined.
Voskresenski\u\i \ established the following property of the Picard module (see \cite[4.48]{Vo-rus})
which is of utmost importance for the whole theory:
\begin{theorem} \label{th:fl-mod}
The module $\Pic \Vb$ is flasque.
\end{theorem}
As mentioned above, this result
gave rise to a purely algebraic way of constructing the flasque resolutions, as well
as other types of resolutions; flasque tori and torsors under such tori were objects
of further thorough investigation \cite{CTS87a}, \cite{CTS87b}.

\begin{remark} \label{rem:bryl}
The geometric method of constructing flasque resolutions described above can be extended
to arbitrary characteristic. This can be done in the most natural way after the gap in Brylinski's proof
\cite{Br} of the existence of a smooth complete model of any torus have been filled in
\cite{CTHS}.
\end{remark}

\begin{remark} \label{rem:fl-mod}
Recently Theorem \ref{th:fl-mod} was extended to the Picard module of a
smooth compactification of an arbitrary connected linear algebraic group \cite{BK04} and,
even more generally, of a homogeneous space of such a group with connected geometric stablizer \cite{CTK06}.
In each case, there is a reduction to the case of tori (although that in \cite{CTK06} is involved enough).
\end{remark}

\subsection{Noether's problem} \label{sec:N}

One of the most striking applications of the birational invariant described above is
a construction of counter-examples to a problem of E.~Noether on rationality of
the field of rational functions invariant under a finite group $G$ of permutations.
Such an example first appeared in a paper by Swan \cite{Sw} where tori were not
mentioned but a resolution of type (\ref{eq:fl}) played a crucial role; at the same time
Voskresenski\u\i \ \cite{Vo69} considered the same resolution to prove that a certain
torus is nonrational. In a later paper \cite{Vo70} he formulated in an explicit way
that the field of invariants under consideration is isomorphic to the function field of
an algebraic torus. This discovery yielded a series of subsequent papers (Endo--Miyata,
Lenstra, and Voskresenski\u\i \ himself) which led to almost complete understanding of the
case where the finite group $G$ acting on the function field is abelian; see \cite[Ch.~VII]{Vo-rus}
for a detailed account. Moreover, the idea of realizing some field of invariants as the function
field of a certain torus proved useful in many other problems of the theory of invariants,
see works of Beneish, Hajja, Kang, Lemire, Lorenz, Saltman, and others;
an extensive bibliography can be found in the monograph \cite{Lo}, see \cite{Ka} and references
therein for some more recent development; note also an alternative approach to Noether's problem
based on cohomological invariants (Serre, in \cite{GMS}).

Note a significant difference in the proofs of nonrationality for the cases $G=\Z_{47}$ (the smallest
counter-example for a cyclic group of prime order) and $G=\Z_8$ (the smallest counter-example for an arbitrary
cyclic group). If $G$ is a cyclic group of order $q=p^n$, and $k$ is a field of characteristic
different from $p$, the corresponding $k$-torus splits over the cyclotomic extension $k(\zeta_q)$.
If $p>2$, the extension $k(\zeta_q)/k$ is cyclic. According to \cite{EM}, if a $k$-torus $T$ splits over
an extension $L$ such that all Sylow subgroups of $\Gal (L/k)$ are cyclic, then the $\Gamma$-module $F$
in the flasque resolution for $T$ is a direct summand of a permutation module. Thus to prove that
it is not a permutation module, one has to use subtle arguments. If $p=2$ and $n\ge 3$, the Galois group
$\Gamma$ of $k(\zeta_q)/k$ may be noncyclic and contain a subgroup $\Gamma'$ such that
$H^1(\Gamma ',F)\ne 0$ which guarantees that $F$ is not a permutation module and hence $T$ is not rational.
This important observation was made in \cite{Vo73}.
Another remark should be done here: for the tori appearing in Noether's problem over a field $k$ with cyclic $G$
such that $\char (k)$ is prime to the exponent of $G$, triviality of the similarity class $[F]$ is necessary and {\it sufficient}
condition for rationality of $T$ (and hence of the corresponding field of invariants). This is an important
instance of the following phenomenon: in a certain class of tori any stably rational torus is rational.
The question whether this principle holds in general is known as Zariski's problem for tori and is left out of the
scope of the present survey.

The group $H^1(\Gamma ,\Pic \Xb )$, where $X$ is a smooth compactification of a $k$-torus $T$, admits another
interpretation: it is isomorphic to $\Br X/\Br k$, where $\Br X$ stands for the Brauer--Grothendieck group of $X$.
This birational invariant, later named the unramified Brauer group, played an important role
in various problems, as well as its generalization for higher unramified cohomology (see \cite{CT96a}, \cite{Sa96}
for details). Let us only note that the main idea here is to avoid explicit construction of a smooth compactification
of an affine variety $V$ under consideration, trying to express $\Br X/\Br k$ in terms of $V$ itself. In the
toric case this corresponds to the formula \cite{CTS87b}
\begin{equation}
\Br X/\Br k = \ker [H^2(\Gamma ,\Th) \to \prod_C H^2(C,\Th )],
\label{eq:sha}
\end{equation}
where the product is taken over all cyclic subgroups of $\Gamma$. Formulas of similar flavour were obtained
for the cases where $V=G$, an arbitrary connected linear algebraic group \cite{CTK98}, \cite{BK00}, and
$V=G/H$, a homogeneous space of a simply connected group $G$ with connected stabilizer $H$ \cite{CTK06};
the latter formula was used in \cite{CTKPR} for proving nonrationality of the field extensions of the
form $k(V)/k(V)^G$, where $k$ is an algebraically closed field of characteristic zero, $G$ is a simple
$k$-group of any type except for $A_n$, $C_n$, $G_2$, and $V$ is either the representation of $G$ on
itself by conjugation or the adjoint representation on its Lie algebra.

It is also interesting to note that the same invariant, the unramified Brauer group of the quotient
space $V/G$, where $G$ is a finite group and $V$ its faithful complex representation, was used by Saltman
\cite{Sa84} to produce the first counter-example to Noether's problem over $\C$. In the same spirit as in formula
(\ref{eq:sha}) above,
this invariant can be expressed solely in terms of $G$: it equals
\begin{equation}
B_0(G):= \ker [H^2(G ,\Q/\Z) \to \prod_A H^2(A,\Q/\Z )],
\label{eq:Bog}
\end{equation}
where the product is taken over all abelian subgroups of $G$ (and,
in fact, may be taken over all bicyclic subgroups of $G$)
\cite{Bo87}. This explicit formula yielded many new counter-examples
(all arising for nilpotent groups $G$, particularly from $p$-groups
of nilpotency class 2). The reader interested in historical
perspective and geometric context is referred to \cite{Sh},
\cite{CTS07}, \cite[6.6, 6.7]{GS}, \cite{Bo07}. We only mention here
some recent work \cite{BMP}, \cite{Ku07} showing that such
counter-examples cannot occur if $G$ is a simple group; this
confirms a conjecture stated in \cite{Bo92}.

\subsection{Generic tori} \label{sec:gen}

Having at our disposal examples of tori with ``good'' and ``bad''
birational properties, it is natural to ask what type of behaviour is
typical. Questions of such ``nonbinary'' type, which do not admit an answer
of the form ``yes-no'', have been considered by many mathematicians, from Poincar\'e to Arnold,
as the most interesting ones. In the toric context, the starting point was the famous
Chevalley--Grothendieck theorem stating that the variety of maximal tori in a connected linear algebraic
group $G$ is rational. If $G$ is defined over an algebraically closed field, its underlying variety is
rational. However, if $k$ is not algebraically closed, $k$-rationality (or nonrationality) of $G$ is a hard
problem. The Chevalley--Grothendieck theorem gives a motivation for studying generic tori in $G$:
if this torus is rational, it gives the $k$-rationality of $G$. The notion of generic torus
can be expressed in Galois-theoretic terms: these are tori whose minimal splitting field has ``maximal possible''
Galois group $\Gamma$ (i.e. $\Gamma$ lies between $W(R)$ and $\Aut (R)$ where $R$ stands for the root system
of $G$). This result, going back to E.~Cartan, was proved in \cite{Vo88}.
In this way, Voskresenski\u\i \
and Klyachko \cite{VoKl} proved the rationality of all adjoint groups of type $A_{2n}$. The rationality was
earlier known for the adjoint groups of type $B_n$, the simply connected groups of type $C_n$, the inner forms
of the adjoint groups of type $A_n$, and (after Theorem \ref{th:two}) for all groups of rank at most two.
It turned out that for the adjoint and simply connected groups of all remaining types the generic torus is not
even stably rational \cite{CK}; in most cases this was proved by computing the birational invariant
$H^1(\Gamma ',F)$ for certain subgroups $\Gamma'\subset\Gamma$. In the case of inner forms of simply connected groups
of type $A_n$, corresponding to generic norm one tori,
this confirms a conjecture by Le Bruyn \cite{LB} (independently proved later in \cite{LL}).
The above theorem  was extended to
groups which are neither simply connected nor adjoint and heavily used in the
classification of linear algebraic groups admitting a rational parameterization of Cayley type \cite{LPR}.

The above mentioned results may give an impression that except for certain types of groups the behaviour
of generic tori is ``bad'' from birational point of view. However, there is also a positive result:
if $T$ is a generic torus in $G$, then $H^1(\Gamma ,F)=0$, This was first proved in \cite{VoKu} for generic tori in the
classical simply connected (type $A_n$) and adjoint groups, and in \cite{Kl89} in the general case. (An independent
proof for the simply connected case was communicated to the author by M.~Borovoi.)
This result has several number-theoretic applications, see Section \ref{sec:arith}.

\begin{remark}
Yet another approach to the notion of generic torus was developed
in \cite{Gr05b} where the author, with an eye towards arithmetic applications,
considered maximal tori in semisimple simply connected groups arising as
the centralizers of regular semisimple stable conjugacy classes.
\end{remark}

\section{Relationship with arithmetic}   \label{sec:arith}

\subsection{Global fields: Hasse principle and weak approximation}

According to a general principle formulated in \cite{Ma},
the influence of (birational) geometry of a variety on its
arithmetic (diophantine) properties may often be revealed
via some algebraic (Galois-cohomological) invariants. In the toric
case such a relationship was discovered by Voskresenski\u\i \  and stated
as the exact sequence (see \cite[6.38]{Vo-rus})
\begin{equation}
0\to A(T) \to H^1(k, \Pic \Vb )\,\tilde{}  \to \Sha (T) \to 0,
\label{eq:a-sh}
\end{equation}
where $k$ is a number field, $T$ is a $k$-torus, $V$ is a
smooth compactification of $T$, $A(T)$ is the defect of weak
approximation, $\Sha (T)$ is the Shafarevich--Tate group,
and $\tilde{}$ stands for the Pontrjagin duality of abelian groups.
The cohomological invariant in the middle, being a purely
algebraic object, governs arithmetic properties of $T$.

On specializing $T$ to be the norm one torus corresponding to
a finite field extension $K/k$, we get a convenient algebraic
condition sufficient for the Hasse norm principle to hold for
$K/k$. In particular, together with results described
in Section \ref{sec:gen}, this shows that the Hasse norm
principle holds ``generically'', i.e. for any field extension
$K/k$ of degree $n$ such that the Galois group of the normal
closure is the symmetric group $S_n$ \cite{VoKu}. (Another
proof was independently found for $n>7$ by Yu.~A.~Drakokhrust.)
Another application was found in \cite{KS}: combining the above mentioned theorem of Klyachko
with the Chevalley--Grothendieck theorem and Hilbert's irreducibility
theorem, one can produce a uniform proof of weak approximation property
for all simply connected, adjoint and absolutely almost simple groups.
(Another uniform proof was found by Harari \cite{Ha} as a consequence
of a stronger result on the uniqueness of the Brauer--Manin obstruction, see
the next paragraph.) Yet another interesting application of sequence (\ref{eq:a-sh})
refers to counting points of bounded height on smooth compactifications of tori
\cite{BT95}, \cite{BT98}: the constant appearing in the asymptotic formula of Peyre \cite{Pe}
must be corrected by a factor equal to the order of $H^1(k, \Pic \Vb )$ and arising on the proof as the
product of the orders of $A(T)$ and $\Sha (T)$ (I thank J.-L.~Colliot-Th\'el\`ene for this remark).

The sequence (\ref{eq:a-sh}) was extended by Sansuc \cite{San} to the case
of arbitrary linear algebraic groups. On identifying the invariant in the middle with
$\Br V/\Br k$, as in Section \ref{sec:N}, one can put this result into more general
context of the so-called Brauer--Manin obstruction to the Hasse principle and weak
approximation (which is thus the only one for principal homogeneous spaces of linear
algebraic groups). This research, started in \cite{Ma}, gave many impressive
results. It is beyond the scope of the present survey.

\begin{remark} \label{rem:strong}
Other types of approximation properties for tori have been considered
in \cite{CTSu} (weaker than weak approximation), \cite{Ra}, \cite{PR}
(strong approximation with respect to certain infinite sets of primes with infinite
complements --- so-called generalized arithmetic progressions). In the latter paper
generic tori described above also played an important role. They were also used
in \cite{CU} in a quite different arithmetic context.
\end{remark}

\subsection{Arithmetic of tori over more general fields}  \label{sec:fields}

Approximation properties and local-global principles were studied
for some function fields. In \cite{An} the exact sequence (\ref{eq:a-sh})
has been extended to the case where the ground field $k$ is pseudoglobal,
i.e. $k$ is a function field in one variable whose field of constants
$\kappa$ is pseudofinite (this means that $\kappa$ has exactly one
extension of degree $n$ for every $n$ and every absolutely irreducible
affine $\kappa$-variety has a $\kappa$-rational point). In \cite{CT96b}
weak approximation and the Hasse principle were established for any torus
defined over $\R (X)$ where $X$ is an irreducible real curve. This allows
one to establish these properties for arbitrary groups over such fields, and,
more generally, over the fields of virtual cohomological
dimension 1 \cite{Sch}.
The same properties for tori defined over some geometric
fields of dimension 2 (such as a function field in two variables over an
algebraically closed field of characteristic zero, or the fraction field of
a two-dimensional. exce;;ent, Henselian, local domain with algebraically closed
residue field, or
the field of Laurent series in one variable over a field of characteristic
zero and cohomological dimension one)
were considered in \cite{CTGP}. Here one
can note an interesting phenomenon: there are counter-examples to
weak approximation but no counter-example to the Hasse principle
is known. One can ask whether there exists some Galois-cohomological
invariant of tori defined over more general fields whose vanishing
would guarantee weak approximation property for the torus under
consideration. Apart from the geometric fields considered in \cite{CTGP},
another interesting case could be $k=\Q_p(X)$, where $X$ is an irreducible
$\Q_p$-curve; here some useful cohomological machinery has been developed in
\cite{SvH}.

\subsection{Integral models and class numbers of tori}
The theory of integral models of tori, started by Raynaud who
constructed an analogue of the N\'eron smooth model \cite{BLR},
has been extensively studied during the past years, and some interesting
applications were found using both N\'eron--Raynaud models and Voskresenski\u\i 's
``standard'' models. The interested reader is referred to the bibliography in
\cite{VKM}. Some more recent works include standard integral models of toric varieties
\cite{KM} and formal models for some classes of tori \cite{DGX}.

Main results on class numbers of algebraic tori are summarized in \cite{Vo-eng}.
One can only add that the toric analogue of Dirichlet's class number formula
established in \cite{Shyr} suggests that a toric analogue of the Brauer--Siegel
theorem may also exist. A conjectural formula of the Brauer--Siegel type for
constant tori defined over a global function field can be found in a recent
paper \cite{KT}.

\section{$R$-equivalence and zero-cycles} \label{sec:R}
$R$-equivalence on the set of rational points of an algebraic variety introduced in \cite{Ma}
turned out to be an extremely powerful birational invariant. Its study in the context
of algebraic groups, initiated in \cite{CTS77}, yielded many striking achievements.
We shall only recall here that the first example of a simply connected group whose
underlying variety is not $k$-rational is a consequence of
an isomorphism, established by  Voskresenski\u\i , between $G(k)/R$, where $G=SL(1,D)$, the group of
norm 1 elements in a division algebra over $k$, with the reduced Whitehead group
$SK_1(D)$; as the latter group may be nonzero because of a theorem by Platonov \cite{Pl}, this
gives the needed nonrationality of $G$. This breakthrough gave rise to dozens of papers
on the topic certainly deserving a separate survey. For the lack of such, the
interested reader is referred to \cite[Ch.~6]{Vo-eng}, \cite[Sections~24--33]{Gi07a}, \cite{Gi07b}.

As to $R$-equivalence on tori, the most intriguing question concerns relationship
between $T(k)/R$ and the group $A_0(X)$ of classes of 0-cycles of degree 0 on
a smooth compactification $X$ of $T$. In a recent paper \cite{Me} Merkurjev
proved that these two abelian groups are isomorphic if $T$ is of dimension 3
(for tori of dimension at most 2 both groups are zero because of their birational
invariance and Theorem \ref{th:two}). For such tori he also obtained a beautiful
formula expressing $T(k)/R$ in ``intrinsic'' terms, which does not require
constructing $X$ as above nor a flasque resolution of $\Th$ as in \cite{CTS77}:
$T(k)/R\cong H^1(k,T^{\circ})/R$, where $T^{\circ}$ denotes the dual torus (i.e.
$\hat{T^{\circ }}=\Hom (\Th ,\Z)$). (In the general case, it is not even known whether
the map $X(k)/R\to A_0(X)$ is surjective, see \cite[\S 4]{CT05} for more details.)
As a consequence, Merkurjev obtained an explicit formula for the Chow group $CH_0(T)$ of
classes of 0-cycles on a torus $T$ of dimension at most 3:
$CH_0(T)\cong T(k)/R\oplus \Z _{i_T}$ where $i_T$ denotes
the greatest common divisor of the degrees of all field extensions $L/k$ such that
the torus $T_L$ is isotropic. The proofs, among other things, use earlier results
\cite{Kl82}, \cite{MP} on the $K$-theory of toric varieties. As mentioned above,
they do not rely on the classification of 3-dimensional tori.

To conclude this section,
one can also add the same references \cite{An}, \cite{CTGP}
as in Section \ref{sec:fields} for $R$-equivalence on tori over more general fields.

\section{Applications in information theory}

\subsection{Primality testing}

One should mention here several recent papers \cite{Gr05a}, \cite{Ki} trying to interpret
in toric terms some known methods for checking whether a given integer $n$ is a prime.
In fact, this approach goes back to a much older paper \cite{CC} where the authors
noticed symmetries in the sequences of Lucas type used in such tests (though algebraic
tori do not  explicitly show up in \cite{CC}).

\subsection{Public-key cryptography}

A new cryptosystem based on the discrete logarithm problem in the group
of rational points of an algebraic torus $T$ defined over a finite field
was recently invented by Rubin and Silverberg \cite{RS03}, \cite{RS04a},
\cite{RS04b}. Since this cryptosystem possesses a compression property,
i.e. allows one to use less memory at the same security level, it drew
serious attention of applied cryptographers and yielded a series of
papers devoted to implementation issues \cite{DGPRSSW}, \cite{DW}, \cite{GV},
\cite{Ko} (in the latter paper another interesting approach is suggested
based on representing a given torus as a quotient of the generalized jacobian
of a singular hyperelliptic curve). In \cite{GPS} the authors propose to use
a similar idea of compression for using tori in an even more recent
cryptographic protocol (so-called pairing-based cryptography). It is interesting
to note that the efficiency (compression factor) of the above mentioned
cryptosystems heavily depends on {\it rationality} of tori under consideration
(more precisely, on an explicit rational parameterization of the underlying
variety). As the tori used by Rubin and Silverberg are known to be stably
rational, the seemingly abstract question on rationality of a given stably rational
torus is moving to the area of applied mathematics. The first challenging problem
here is to obtain an explicit rational parameterization of the 8-dimensional
torus $T_{30}$, defined over a finite field $k$ and splitting over its cyclic extension $L$
of degree 30, whose character module $\Th _{30}$ is isomorphic to $\Z[\zeta_{30}]$,
where $\Z [\zeta_{30}]$ stands for a primitive 30th root of unity. (Here is an alternative
description of $T_{30}$: it is a maximal torus in $E_8$ such that $\Gal (L/k)$ acts on $\Th _{30}$
as the Coxeter element of $W(E_8)$; this can be checked by a direct computation or using \cite{BF}.)

This is a
particular case of a problem posed by Voskresenski\u\i \ \cite[Problem 5.12]{Vo-rus}
30 years ago. Let us hope that we will not have to wait another 30 years for answering
this question on a degree 30 extension.

\medskip

\noindent {\it Acknowledgements}. The author's research was
supported in part by the Minerva Foundation through the Emmy Noether
Research Institute of Mathematics and by a grant from the Ministry
of Science, Culture and Sport, Israel, and the Russian Foundation
for Basic Research, the Russian Federation. This paper was written
during the visit to the MPIM (Bonn) in August--September 2007. The
support of these institutions is highly appreciated. I thank
J.-L.~Colliot-Th\'el\`ene for many helpful remarks.

\end{document}